%% file: vanish.tex
\documentclass{amsart}
 
 \usepackage{amscd}
 
 \input skip-thm.tex

\input qformsl-def.tex

 \theoremstyle{plain}
\newtheorem{specthmp}{Theorem}

 \numberwithin{equation}{section}

 \DeclareMathOperator{\trace}{trace}
 \DeclareMathOperator{\rad}{rad}

 \newcommand{\roots}{R}

\DeclareMathOperator{\PSp}{PSp}
\DeclareMathOperator{\PSO}{PSO}
\DeclareMathOperator{\HSpin}{HSpin}

\newcommand{\RG}{RG}

\newcommand{\co}{{\vee}}
\newcommand{\ach}{\alpha^\co}
\newcommand{\mch}{m^\co}
\newcommand{\hch}{h^\co}
\newcommand{\dch}{\delta^\co}

\newcommand{\coroots}{\roots^\co}
\newcommand{\Pch}{P^\co}
\newcommand{\Qch}{Q^\co}

\newcommand{\at}{\tilde{\alpha}}
\newcommand{\atch}{\tilde{\alpha}^\co}

\newcommand{\g}{\Lie(G)}

\renewcommand{\c}{\mathfrak{c}}

\newcommand{\Gt}{\widetilde{G}}
\newcommand{\Gb}{\bar{G}}
\newcommand{\Tb}{\bar{T}}
\newcommand{\Tt}{\widetilde{T}}

\newcommand{\bt}{\widetilde{b}}
\newcommand{\qt}{\widetilde{q}}

\renewcommand{\sl}{\mathfrak{sl}}

\newcommand{\drho}{{\mathrm{d}\rho}}
\newcommand{\df}{{\mathrm{d}f}}
\newcommand{\dpi}{{\mathrm{d}\pi}}

\newcommand{\gl}{\mathfrak{gl}}
 
\begin{document}
 \title[Vanishing of trace forms]{Vanishing of trace forms in low characteristics}
 
 \dedicatory{With an appendix by Alexander Premet}
  
\author{Skip Garibaldi}
\address{(Garibaldi) Department of Mathematics \& Computer Science, Emory University, Atlanta, GA 30322, USA}
\email{skip@member.ams.org}
\urladdr{http://www.mathcs.emory.edu/{\textasciitilde}skip/}

\address{(Premet) School of Mathematics, The University of Manchester, Oxford Rd, M13 9PL, UK.}
\email{sashap@maths.man.ac.uk}




\setlength{\unitlength}{.75cm}

\begin{abstract}
Every finite-dimensional representation of an algebraic group $G$ gives a trace symmetric bilinear form on the Lie algebra of $G$.  We give criteria in terms of root system data for the existence of a representation such that this form is nonzero or nondegenerate.
As a corollary, we show that a Lie algebra of type $E_8$ over a field of characteristic 5 does not have a so-called ``quotient trace form", answering a question posed in the 1960s.
\end{abstract}

\maketitle

Let $G$ be an algebraic group over a field $F$, acting on a finite-dimensional vector space $V$ via a homomorphism $\rho \!: G \ra GL(V)$.  The differential $\drho$ of $\rho$ maps the Lie algebra $\g$ of $G$ into $\mathfrak{gl}(V)$, and we put $\Tr_\rho$ for the symmetric bilinear form
\[ 
\Tr_\rho(x, y) := \trace (\drho(x) \, \drho(y)) \quad \text{for $x, y \in \g$.}
\] 
We call $\Tr_\rho$ a \emph{trace form of $G$}.   Such forms appear, for example, in the hypotheses for the Jacobson-Morozov Theorem \cite[5.3.1]{Carter:big}, in Richardson's proof that there are finitely many conjugacy classes of nilpotent elements in the Lie algebra of a semisimple algebraic group as in \cite[\S2]{Jantzen:nil} or \cite[\S\S3.8, 3.9]{Hum:cc}, and in the ``explicit" construction of a Springer isomorphism in \cite[\S9.3]{BR}.
We prove:

\begin{specthm} \label{MT}
Let $G$ be a split and almost simple linear algebraic group over a field $F$.  Then:
\begin{enumerate}
\item There is a representation $\rho$ of $G$ with $\Tr_\rho$ nondegenerate if and only if the characteristic of $F$ is very good for $G$. \label{MT.nd}
\item There is a representation $\rho$ of $G$ with $\Tr_\rho$ nonzero if and only if the characteristic of $F$ is as indicated in Table \ref{bad.primes}. \label{MT.nz}
\end{enumerate}
\end{specthm}

\begin{table}[t]
\begin{tabular}{p{1.9in}||p{1.2in}|p{1.545in}} 
{\ \newline\centerline{$G$}}&Every $\rho$ has $\Tr_\rho$ degenerate if $\chr F$ is&Every $\rho$ has $\Tr_\rho$ zero if $\chr F$ is\\ \hline\hline
$\SL_n / \mmu{m}$ with $m$ odd&a divisor of $n$&a divisor of $\gcd(m, n/m)$\\
$\SL_n / \mmu{m}$ with $m$ even&a divisor of $n$&a divisor of $2 \gcd(m, n/m)$\\
$\Sp_{2n}$ & 2&none \\
$\left. \parbox{1.8in}{$\SO_n$, $\Spin_n$, $\PSO_n$ with $n \ge 7$\\
$\HSpin_{4n}$ for $n \ge 3$, $\PSp_{2n}$} \right\}$&2&2 \\
$E_6$ adjoint, $G_2$&2 or 3&2 \\
$E_6$ simply connected,
$E_7$, $F_4$&2 or 3&2 or 3 \\
$E_8$&2, 3, or 5&2, 3, or 5
\end{tabular}
\medskip
\caption{Primes where $\Tr_\rho$ is degenerate or zero for every $\rho$.  The middle column lists the primes that are not very good for $G$.  For simply connected $G$, the right column lists the torsion primes for $G$ as defined in, e.g., \cite[1.13]{St:tor}.} \label{bad.primes}
\end{table}

A weaker version (``up to isogeny") of the ``if" direction of part \eqref{MT.nd} is standard; see e.g.\ \cite[I.5.3]{SpSt} or \cite[1.16]{Carter:big}.  
After this paper was released as a preprint, I learned that Alexander Premet had previously proved the ``only if" direction of part \eqref{MT.nd} for groups not of type $A_n$, cf.~\cite[p.~80]{Premet:KW}, but his proof has not been published.  
Here both directions of part \eqref{MT.nd} are deduced from part \eqref{MT.nz}.
The crux of the proof of part \eqref{MT.nz} is a formula for the trace form $\Tr_\rho$, given in Proposition \ref{divides} below.  

We remark that the characteristics in part \eqref{MT.nd} of the theorem depend only on the isogeny class of $G$, whereas the 
characteristics in part \eqref{MT.nz} of the theorem are more sensitive.  For example, $\Sp_{2n}$ has 
a representation with nonzero trace form over every field,
whereas its quotient $\PSp_{2n}$ has such only in characteristic different from 2; this is no surprise because $\PSp_{2n}$ has ``fewer" representations than $\Sp_{2n}$. 
But the opposite phenomenon also occurs: a simply connected group of type $E_6$ has a representation with nonzero trace form only in characteristics $\ne 2, 3$, whereas its quotient the adjoint $E_6$ has one over every field of characteristic $\ne 2$.  This ``opposite phenomenon" is related to the number $E(G)$ defined in \S\ref{E.sec} below.

\smallskip
For $G$ of type $E_8$, we can strengthen Theorem \ref{MT}.  Given a representation $\psi$ of the Lie algebra $\Lie(G)$ of $G$, one can define a trace form $\Tr_\psi$ on $\Lie(G)$ by setting $(x, y) \mapsto \trace(\psi(x)\,\psi(y))$.  We prove:

\begin{specthm} \label{alg}
If $F$ has characteristic $2$, $3$, or $5$ and $G$ is of type $E_8$, then the trace form of every representation of $\Lie(G)$ is zero.
\end{specthm}

This is a strengthening of Theorem \ref{MT} because ``many" representations $\psi$ of $\Lie(G)$ are not differentials of representations of $G$.  The proof of Theorem \ref{alg} is given in \ref{alg.pf}; it amounts to a combination of Theorem \ref{MT} and a result generously provided by Premet, presented in an appendix.  The converse of Theorem \ref{alg} is of course true; in characteristic $\ne 2, 3, 5$, the Killing form is nondegenerate.

In characteristic 5,  Theorem \ref{alg} easily gives an apparently stronger statement, namely that $\Lie(G)$ has no ``quotient trace form", see Cor.~\ref{E8.cor}.  This answers a question posed in the early 1960s, see e.g.\ \cite[p.~554]{Block:trace}, \cite[p.~543]{BlockZ}, or \cite[p.~48]{Sel:mod}. 

From the point of view of Lie algebras, this paper addresses the existence of restricted representations with nonzero or nondegenerate trace forms on Lie algebras of almost simple algebraic groups.  These algebras are approximately the simple Lie algebras ``of classical type".  For fields of characteristic $\ge 5$ and simple Lie algebras of other types  (necessarily Cartan or Melikian by Block-Premet-Strade-Wilson, see \cite{Strade} or \cite{Mathieu}), every representation has zero trace form by \cite[Cor.~3.1]{Block:trace}.



\subsection*{Notation}
All algebraic groups discussed here are linear.  Such a group $G$ over a field $F$ is \emph{almost simple} if it is semisimple and has no proper connected, closed, normal subgroups defined over $F$.  In case $F$ is separably closed, the almost simple algebraic groups are the semisimple groups whose Dynkin diagrams are connected.

$\PSO_n$ denotes the adjoint group of the (split) special orthogonal group $\SO_n$; when $n$ is odd it is the same as $\SO_n$.  Similarly, $\PSp_{2n}$ is the adjoint group of type $C_n$; it can be viewed as $\Sp_{2n} / \mmu2$.  
The groups $\SO_n$, $\Spin_n$, and $\PSO_n$ for $n = 3, 5$, and 6 are isogenous to $\SL_2$, $\Sp_4$, or $\SL_4$ and appear in Table \ref{bad.primes} in that alternative form.
For $n \ge 3$, we write $\HSpin_{4n}$ for the nontrivial quotient of $\Spin_{4n}$ that is neither $\SO_{4n}$ nor adjoint.  

\medskip
{\noindent{\small{\textbf{Acknowledgments.}} It is a pleasure to thank Jean-Pierre Serre for numerous helpful discussions and for pointing out the utility of the form $\bt$ (from \cite{GrossNebe}, see \ref{b.def}) early in this project.  I thank also George McNinch and Burt Totaro for their comments, the NSF for its support, and the Institut des Hautes Etudes Scientifiques for its hospitality.}}

\section{The number $N(G)$ and the Dynkin index} \label{N.sec}

\begin{borel*} \label{N.def}
Fix a simple root system $\roots$.  We write $P$ for its weight lattice and $\qform{\ ,\ }$ for the canonical pairing between $P$ and its dual.
Fix a long root $\alpha \in \roots$ and write $\ach$ for the associated coroot.  
For each subset $X$ of $P$ that is invariant under the Weyl group, we put:
\[
N(X) := \frac12 \sum_{x \in X} \qform{x, \ach}^2 \quad \in  \Z\!\left[ \textstyle\frac12 \right].
\]
The number \emph{$N(X)$ does not depend on the choice of $\alpha$} because the long roots are conjugate under the Weyl group.

Furthermore, \emph{$N(X)$ is an integer.}  To see this, note that the reflection $s$ in the hyperplane orthogonal to $\alpha$ satisfies
$\qform{s x, \ach} = \qform{x, s \ach} = -\qform{x, \ach}$,
so in the definition of $N(X)$, the sum can be taken to run over those $x$ satisfying $x \ne sx$.  For such $x$, we have
$\qform{x, \ach}^2 + \qform{sx, \ach}^2 = 2 \qform{x, \ach}^2$,
proving the claim.
\end{borel*}

\begin{eg} \label{Ad.eg}
The computations in \cite[pp.~180, 181]{SpSt} show that $N(R) = 2 \hch$, where $\hch$ denotes the \emph{dual Coxeter number} of $\roots$, which is defined as follows.  Fix a set of simple roots $\D$ of $R$.  Write $\at$ for the highest root; the corresponding coroot $\atch$ is
\[
\atch = \sum_{\delta \in \D} \mch_\delta \dch
\]
for some natural numbers $\mch_\delta$.  The dual Coxeter number $\hch$ is defined by
\[
\hch := 1 + \sum \mch_\delta.
\]
In case all the roots of $R$ have the same length, it is the (usual) Coxeter number $h$ and is given in the tables in \cite{BLie}.

Suppose that there are two different root lengths in $R$; we write $L$ for the set of long roots and $S$ for the set of short roots.  The arguments in \cite{SpSt} are easily adapted to show that
\[
N(L) = 2 \left[ 1 + \sum_{\delta \in \D \cap L} \mch_\delta \right] \quad \text{and} \quad N(S) = 2 \sum_{\delta \in \D \cap S} \mch_\delta.
\]
We obtain the following numbers:
\[
\begin{tabular}{c|cccc} 
type of $R$&$h$&$\hch$&$N(L)$&$N(S)$ \\ \hline
$B_n$ ($n \ge 2$)&$2n$&$2n-1$&$4(n-1)$&2 \\
$C_n$ ($n \ge 2$)&$2n$&$n+1$&4&$2(n-1)$ \\
$G_2$&6&4&6&2 \\
$F_4$&12&9&12&6
\end{tabular}
\]
\end{eg}

\begin{defn} \label{Nrho.def}
Fix a split almost simple linear algebraic group $G$ over $F$.  Fix also a pinning of $G$ with respect to some maximal torus $T$; this includes a root system $\roots$ and a set of simple roots $\D$ of $G$ with respect to $T$.
For a representation $\rho$ of $G$ over $F$, one defines
\[
N(\rho) := \sum_{\text{dominant weights $\la$}} \left( \parbox{1in}{multiplicity of $\la$ as a weight of $\rho$} \right) \cdot N(W\la) \quad \in \Z.
\]
For example, the adjoint representation $\Ad$ has $N(\Ad) = 2\hch$ by Example \ref{Ad.eg}.    The number $N(\rho)$ is the \emph{Dynkin index} of the representation $\rho$ defined in \cite[p.~130]{Dynk:ssub} and studied in \cite{MG}.
The Dynkin index of the fundamental irreducible representations of $G$ (over $\C$) are listed in \cite[Prop.~2.6]{LaszloSorger} or \cite[pp.~36--44]{McKPR}, correcting some small errors in Dynkin's calculations.

We put:
\begin{equation} \label{NG.def}
N(G) := \gcd N(\rho),
\end{equation}
where the $\gcd$ runs over the representations of $G$ defined over $F$.
Because the map $\rho \mapsto N(\rho)$ depends only on the weights of $\rho$ with multiplicity, it is compatible with short exact sequences 
\begin{equation} \label{rho.seq}
\begin{CD}
0 @>>> \rho' @>>> \rho @>>> \rho/\rho' @>>> 0
\end{CD}
\end{equation}
in the sense that
\[
N(\rho) = N(\rho') + N(\rho/\rho').
\]
Writing $\RG$ for the representation ring of $G$, we obtain a homomorphism of abelian groups $N \!: \RG \ra \Z$ with image $N(G) \cdot \Z$.

In the definition of $N(G)$, it suffices to let the gcd run over generators of $\RG$, e.g., the irreducible representations of $G$.  For an irreducible representation $\rho$, the highest weight $\la$ has multiplicity 1 and all the other weights of $\rho$ are lower in the partial ordering.  Inducting on the partial ordering, we find:
\[
N(G) = \gcd \{ N(W\la) \mid \la \in T^* \}.
\]
In particular, $N(G)$ depends only on the root system $R$ and the lattice $T^*$, and not on the field $F$.  
\end{defn}

\begin{eg} \label{dynk.def}
When $G$ is simply connected, the number $N(G)$ is known as the \emph{Dynkin index of $G$} and its value is listed in, e.g., \cite{LaszloSorger}.  Examining the list of values, one finds that the primes dividing $N(G)$ (for $G$ simply connected) are the torsion primes of $G$.
\end{eg}

\begin{eg} \label{SO.N}
Write $\Spin_n$ and $\SO_n$ for the spin and special orthogonal groups of an $n$-dimensional nondegenerate quadratic form of maximal Witt index.  For $n \ge 7$, these groups are split and almost simple of type $B_\ell$ (with $\ell \ge 3$) or $D_\ell$ (with $\ell \ge 4$).  The Dynkin index $N(\Spin_n)$ is 2; it obviously divides $N(\SO_n)$.  On the other hand, the natural $n$-dimensional representation $\rho$ of $\SO_n$ has $N(\rho) = 2$, so $N(\SO_n) = 2$.
\end{eg}

\begin{eg} \label{PSp.N}
We claim that 
\[
N(\PSp_{2n}) = \begin{cases}
2&\text{if $n$ is even}\\
4&\text{if $n$ is odd}
\end{cases}
\]
for $n \ge 2$.
The number $N(\PSp_{2n})$ divides 4 and $2(n-1)$ by Example \ref{Ad.eg}.  Further, $N(\PSp_{2n})$ is even by \cite[14.2]{MG}.  This shows that $N(\PSp_{2n})$ is 2 or 4, and is 2 in case $n$ is even.

Suppose that $n$ is odd.  We must show that $N(W\la)$ is divisible by 4 for every element $\la$ of the root lattice of $\PSp_{2n}$.  We use the same notation as \cite[\S14]{MG} for the weights of $\PSp_{2n}$: they are a sum $\sum_{i=1}^n x_i e_i$ such that $\sum x_i$ is even.  The Weyl group $W$ is a semidirect product of $(\Zm2)^n$ (acting by flipping the signs of the $e_i$) and the symmetric group on $n$ letters (acting by permuting the $e_i$).  Taking $X$ for the $(\Zm2)^n$-orbit of $\sum x_i e_i$, we have
\begin{equation} \label{M.144}
\frac12 \sum_{x \in X} \qform{ \sum\nolimits_i x_i e_i,  (2 e_n)^\vee}^2 = 2^{r-1} x_n^2
\end{equation}
where $r$ denotes the number of nonzero $x_i$'s, cf.~\cite[pf.~of Lemma 14.2]{MG}.  If $r = 1$, then the unique nonzero $x_i$ is even, and we find that for $r \ne 2$, the sum---hence also $N(W \sum x_i e_i)$---is divisible by 4.  Suppose that $x_1, x_2$ are the only nonzero $x_i$'s; then by \eqref{M.144} we have:
\[
N(W(x_1 e_1 + x_2 e_2)) = \begin{cases}
2(n-1)(x_1^2 + x_2^2) & \text{if $x_1 \ne \pm x_2$} \\
2(n-1)x_1^2&\text{if $x_1 = \pm x_2$.}
\end{cases}
\]
As $n$ is odd,
$N(W(x_1 e_1 + x_2 e_2))$ is divisible by 4, which completes the proof of the claim.  
\end{eg}

\begin{eg} \label{E7.N}
\emph{For $G$ adjoint of type $E_7$, we have $N(G) = 12$.}  
To see this, we note that $N(G)$ is divisible by $N(\Gt)$, where $\Gt$ is the universal covering of $G$, i.e., 12 divides $N(G)$.  Also, $N(G)$ divides $2 \hch = 36$ by Example \ref{Ad.eg}.  
For the minuscule representation $\rho$ of $\Gt$, we have $\dim \rho = 56$ and $N(\rho) = 12$.  The representation $\rho^{\ot 2}$ of $\Gt$ factors through $G$ and by 
the ``derivation formula" 
\[
N(\rho_1 \ot \rho_2) = (\dim \rho_1) \cdot N(\rho_2) + (\dim \rho_2) \cdot N(\rho_1)
\]
(see e.g.\ \cite[p.~122]{MG}) we have 
\[
N(\rho^{\ot 2}) = 2 (\dim \rho) N(\rho) = 2^6 \cdot 3 \cdot 7.
\]
It follows that $N(G)$ equals 12, as claimed.
\end{eg}

\section{The Lie algebra of $G$} \label{Lie.sec}
 
\begin{borel*} \label{Lie.def}
Let $G$ be a split almost simple algebraic group over $F$; we fix a pinning for it.  If $G_\Z$ is a split group over $\Z$ with the same root datum as $G$, the pinning identifies $G$ with the group obtained from $G_\Z$ by the base change $\Z \ra F$ and the maximal torus $T$ in $G$ (from the pinning) with the base change of a maximal torus $T_\Z$ in $G_\Z$.  We have a root space decomposition of the Lie algebra of $G_\Z$:
\begin{equation} \label{Cartan}
\Lie(G_\Z) = \Lie(T_\Z) \oplus \bigoplus\nolimits_{\alpha \in R} \Z x_\alpha
\end{equation}
and 
\begin{equation} \label{Lie.T}
\Lie(T_\Z) =  \{ h \in \Lie(T_\C) \mid \text{$\mu(h) \in \Z$ for all $\mu \in T^*$} \},
\end{equation}
see \cite[p.~64]{St}.  Because $\Lie(G_\Z)$ is a free $\Z$-module, the Lie algebra $\Lie(G)$ of $G$ is naturally identified 
with $\Lie(G_\Z) \ot_\Z F$, and similarly for $\Lie(T)$, see \cite[II.4.4.8]{DG}.
\end{borel*}

\begin{borel*}
Write $\Gt$ for the universal covering of $G$; we use the obvious analogues of the notations in \ref{Lie.def} for $\Gt$.  The group $G$ acts on $\Gt$ by conjugation, hence also on $\Lie(\Gt)$.  
If the kernel of the map $\Gt \ra G$ is \'etale, then the representation $\Lie(\Gt)$ is equivalent to the adjoint representation on $\Lie(G)$. 
 In any case, the natural map $\Lie(\Gt) \ra \Lie(G)$ is an isomorphism on the $F$-span of the $x_\alpha$'s.
\end{borel*}

\begin{borel*} \label{Lie.sc}
We claim that \emph{$\Lie(\Gt)$ is a Weyl module for $G$} in the sense of \cite[p.~183]{Jantzen}, i.e., its character is given by Weyl's formula and it is generated as a $G$-module by a highest weight vector.  The first condition holds by \eqref{Cartan}, so it suffices to check the second.  

To check that the submodule $Gx_{\at}$ generated by the highest weight vector $x_{\at}$ is all of $\Lie(\Gt)$, one quickly reduces to checking that $Gx_{\at}$ contains $\Lie(\Tt_\Z)$.  Equation \eqref{Lie.T} gives a natural isomorphism $\Z[\coroots] \iso \Lie(\Tt_\Z)$ where $\Tt_\Z$ is the maximal torus in $\Gt_\Z$ mapping onto $T_\Z$.  We write (as is usual) $h_\alpha$ for the image of $\ach$ under this map.  As $[x_\alpha, x_{-\alpha}] = h_\alpha$, the claim is proved.
\end{borel*}

\begin{borel*} \label{Lie.simple}
See \cite{Hiss}, \cite[esp.~Cor.~2.7a]{Hogeweij}, or \cite[\S2]{St:aut} for descriptions of the composition series of $\Lie(\Gt)$.  They immediately give: \emph{If the characteristic of $F$ is very good for $G$, then $\Lie(G)$ is a simple Lie algebra.  If additionally $F$ is infinite then $\Lie(G)$ is an irreducible $G(F)$-module.}
\end{borel*}

\section{The number $E(G)$}\label{E.sec}

\begin{defn} \label{b.def}
Maintain the notation of the preceding section.  The Killing form on $\Lie(\Gt_\Z)$ is divisible by $2\hch$ \cite{GrossNebe} and dividing by $2\hch$ gives an indivisible even symmetric bilinear form $\bt$ on $\Lie(\Gt_\Z)$ such that
\begin{equation} \label{xlong}
\bt(h_\alpha, h_\alpha) = 2 \eand  \bt(x_\alpha, x_{-\alpha}) = 1
\end{equation}
for long roots $\alpha$, see \cite[p.~181]{SpSt} or \cite[Lemma VIII.2.4.3]{BLie}.  For a short root $\beta$, we have:
$\bt(h_\beta, h_\beta) = 2c$ and $\bt(x_\beta, x_{-\beta}) = c$,
where $c$ is the square-length ratio of $\alpha$ to $\beta$.  For example, $G = \SL_n$ has Lie algebra the trace zero $n$-by-$n$ matrices, and the form $\bt$ is the usual trace bilinear form $(x,y) \mapsto \trace(xy)$, cf.~\cite[Exercise VIII.13.12]{BLie}.

The natural map $\Lie(\Gt_\Z) \ra \Lie(G_\Z)$ is an inclusion and extending scalars to $\Q$ gives an isomorphism.  Therefore, $\bt$ gives a rational-valued symmetric bilinear form on $\Lie(G_\Z)$.  We define $E(G)$ to be the smallest positive rational number such that $E(G)\cdot \bt$ is integer-valued on $\Lie(G_\Z)$; we write $b$ for this form.  Note that $E(G)$ is an integer by \eqref{xlong}.

Clearly, $E(G)$ depends only on the root system of $G$ and the character lattice $T^*$ viewed as a sublattice of the weight lattice, and not on the field $F$.
\end{defn}

\begin{borel*} \label{adj.E}
Write $\Gb$ for the adjoint group of $G$; we use the obvious analogues of the notations in \ref{Lie.def} for $\Gb$.  We have a commutative diagram
\[
\begin{CD}
\Qch @>{\sim}>> \Lie(\Tt_\Z) \\
@VVV @VVV \\
\Pch @>{\sim}>> \Lie(\Tb_\Z)
\end{CD}
\]
where $\Qch$ and $\Pch$ are the root and weight lattices of the dual root system.
The form $\bt$ restricts to be an inner product on $\Qch$ such that the square-length of a short coroot $\ach$ is 2.  This inner product extends to a rational-valued inner product on $\Pch$, and $E(\Gb)$ is the smallest positive integer such that $E(\Gb)\cdot \bt$ is integer-valued on $\Pch$.
\end{borel*}

\begin{eg} \label{E.C}
Consider the case where $G$ is $\PSp_{2n}$ for some $n \ge 2$, i.e., adjoint of type $C_n$.  In the notation of the tables in \cite{BLie}, the form $\bt$ is twice the usual scalar product, i.e., $\bt(e_i, e_j) = 2\delta_{ij}$ (Kronecker delta).  The fundamental weight $\omega_n$ has $\bt(\omega_n, \omega_n) = n/2$.  Checking $\bt(\omega_i, \omega_j)$ for all $i, j$, shows that $E(\Gb)$ is 1 if $n$ is even and 2 if $n$ is odd.
\end{eg}

\begin{eg} \label{E.sl}
Suppose that all the roots of $G$ have the same length, so that we may identify the root system $R$ with its dual and normalize lengths so that $\qform{\, , \, }$ is symmetric and equals $\bt$ on $\Qch$.

(1): \emph{$E(\Gb)$ is the exponent of $P/Q$}, the weight lattice modulo the root lattice.  Indeed, the isomorphism between $P$ and $\Lie(\Tb_\Z)$ shows that $E(\Gb)$ is the smallest natural number such that $E(\Gb) \cdot \qform{\, , \, }$ is integer-valued on $P \times P$, equivalently, the smallest natural number $e$ such that $e P$ is contained in $Q$; this is the exponent of $P/Q$.

(2): The bilinear form
\[
\bt \!: \Lie(\Gt_\Z) \times \Lie(\Gb_\Z) \ra \Q
\]
has image $\Z$ and identifies $\Lie(\Gb_\Z)$ with $\Hom_{\Z}(\Lie(\Gt_\Z), \Z)$.  (On the span of the $x_\alpha$'s, this is clear from \eqref{xlong}.  On the Cartan subalgebras, it amounts to the statement that $\qform{\, , \,}$ identifies $P$ with $\Hom(Q,\Z)$.)
It follows that $\Lie(\Gb)$, as a $G$-module, is the dual of $\Lie(\Gt)$, i.e., $\Lie(\Gb)$ is the module denoted by $H^0(\at)$ in \cite{Jantzen}.
\end{eg}

\begin{eg} \label{SO.E}
For $n = 3$ or $n \ge 5$, we claim that $E(\SO_n) = 1$.  

For $n$ odd, $\SO_n$ is adjoint of type $B_\ell$ for $\ell = (n-1)/2$, and we compute as in \ref{adj.E} and Example \ref{E.C}.  The dual root system is of type $C_\ell$, and the form $\bt$ is the usual scalar product, i.e., $\bt(e_i, e_j) = \delta_{ij}$.  The fundamental weight $\omega_i$ is $e_1 + e_2 + \cdots + e_i$, so $E(\SO_{2\ell+1}) = 1$.

For $n$ even, $\SO_n$ has type $D_\ell$ for $\ell = n/2$.  The character group $T^*$ of a maximal torus in $\SO_n$ consists of the weights whose restriction to the center of $\Spin_n$ is 0 or agrees with the vector representation, i.e., the weights $\sum c_i \omega_i$ such that $c_{\ell-1} + c_\ell$ is even.  It follows that that the cocharacter lattice $T_*$ is generated by the (co)root lattice and 
\[
\omega_1 = \alpha_1 + \alpha_2 + \cdots + \alpha_{\ell - 2} + \frac12 (\alpha_{\ell - 1} + \alpha_\ell).
\]
We have:
\[
\bt(\omega_1, \omega_1) = 
\qform{\omega_1, \omega_1} = 1,
\]
so the form $\bt$ is integer-valued on $T_*$ and $E(\SO_{2\ell}) = 1$.
\end{eg}

\begin{eg} \label{HSpin.E}
Let $\HSpin_{4n}$ denote the image of $\Spin_{4n}$ under the irreducible representation with highest weight $\omega_\ell$ for $\ell := 2n$; it is a ``half-spin" group.  The character lattice $T^*$ consists of weights $\sum c_i \omega_i$ such that $c_1 + c_{\ell - 1}$ is even.  The lattice generated by $Q$ and 
\[
\omega_1 + \omega_{\ell - 1} = \frac12 (3\alpha_1 + 4\alpha_2 + \cdots + \ell \alpha_{\ell - 2}) + \frac{\ell+2}4 \alpha_{\ell - 1} + \frac{\ell}4 \alpha_\ell
\]
contains $Q$ with index 2 and is contained in $T_*$, hence equals $T_*$.  As
\[
\bt(\omega_1 + \omega_{\ell - 1}, \omega_1 + \omega_{\ell - 1}) =
\qform{\omega_1 + \omega_{\ell - 1}, \omega_1 + \omega_{\ell - 1}} = \frac32 + \frac{\ell+2}4 = \frac{n}2 + 2,
\]
we conclude that:
\[
E(\HSpin_{4n}) = \begin{cases}
1 & \text{if $n$ is even} \\
2 & \text{if $n$ is odd.}
\end{cases}
\]
\end{eg}

\section{Formula for the trace} \label{crux.sec}

The integer-valued symmetric bilinear form $b$ on $\Lie(G_\Z)$ defined in \ref{b.def} gives by scalar extension a symmetric bilinear form on $\Lie(G)$ which we denote by $b_F$.

\begin{prop} \label{divides}
Let $\rho$ be a representation of a split and almost simple algebraic group $G$ over $F$.  Then:
\begin{enumerate}
\item $E(G)$ divides $N(\rho)$.
\item $\Tr_\rho = \displaystyle\frac{N(\rho)}{E(G)} b_F$.
\end{enumerate}
\end{prop}

\begin{proof}
We first suppose that $F$ is the complex numbers.  Write $\pi \!: \Gt \ra G$ for the universal covering of $G$ as in \S\ref{Lie.sec}.  We compute $\Tr_{\rho\pi}$.  If we decompose the representation $\rho$ with respect to the action of $\Tt$ and write $V_\mu$ for the eigenspace relative to the weight $\mu$, then $h_\alpha$ acts on $V_\mu$ by scalar multiplication by $\qform{\mu, \ach}$, hence $\Tr_\rho(h_\alpha, h_\alpha) = \sum \dim(V_\mu) \qform{\mu, \ach}^2$.  By putting together the $\mu$ in an orbit $W\la$ (where $\la$ is dominant) and taking $\alpha$ to be a long root, one gets:
\begin{equation} \label{tr.N}
\Tr_{\rho\pi}(h_\alpha, h_\alpha) = 2\, N(\rho).
\end{equation}
The representation $\Lie(\Gt_\Z) \ot \C$ is irreducible and has a nondegenerate $\Gt_{\C}$-invariant symmetric bilinear form, so by Schur's Lemma we have:
\[
\Hom_{\Gt_{\C}}(\Lie(\Gt) \ot \C, (\Lie(\Gt)^*) \ot \C) = \C.
\]
In particular, $\Tr_{\rho\pi}$ equals $z\,\bt$ for some complex number $z$ and
\[
2 N(\rho) = \Tr_{\rho\pi}(h_\alpha, h_\alpha) = z \, \bt(h_\alpha, h_\alpha) = 2z.
\]
Hence $\Tr_{\rho\pi} = N(\rho) \, \bt$.    (This argument can be viewed as restating pp.~130--131 of \cite{Dynk:ssub}.)  The Lie algebra $\Lie(G_\Z) \ot \C$ is naturally identified with $\Lie(\Gt_\Z) \ot \C$, and (2) follows from the equation $E(G) \, \bt = b$ in the case $F = \C$.

Now allow $F$ to be arbitrary but suppose that $\rho$ is a Weyl module.  There is a $\Z$-form $\rho_\Z$ of $\rho$, and the form $\Tr_{\rho_\Z}$ is the restriction of $\Tr_{\rho_\C}$ on $\Lie(G_\C)$ to $\Lie(G_\Z)$.  Because (2) holds over the complex numbers, it holds over the integers, and by scalar extension it holds over the field $F$ as well.  Clearly the form $\Tr_{\rho_\Z}$ is integer valued; as $b$ is indivisible, it follows that $E(G)$ divides $N(\rho)$. 

\smallskip
We now treat the case of an arbitrary representation $\rho$.
The number $N(\rho)$ depends only on the class of $\rho$ in the representation ring $\RG$.  As the Weyl modules generate $\RG$ as an abelian group and $E(G)$ divides $N(\psi)$ for every Weyl module $\psi$, (1) follows.

For (2), we note that the map $\rho \mapsto \Tr_\rho - (N(\rho)/E(G)) \, b_F$ is compatible with exact sequences like \eqref{rho.seq} in the sense that $\Tr_\rho = \Tr_{\rho'} + \Tr_{\rho/\rho'}$. We obtain a homomorphism of abelian groups
\[
\RG \ra \fbox{\parbox{2.1in}{symmetric bilinear forms on $\Lie(G)$}}
\]
that vanishes on the Weyl modules, hence is zero.
\end{proof}


\begin{borel*}
Because $b$ is indivisible (as a form over $\Z$), the form $b_F$ is not zero.  Proposition \ref{divides}(2) immediately gives:
\begin{equation}  \label{single}
\parbox{4.5in}{\emph{Let $\rho$ be a representation of a split and almost simple algebraic group $G$ over $F$.  Then $\Tr_\rho$ is zero if and only if the characteristic of $F$ divides $N(\rho)/E(G)$.}}
\end{equation}
Furthermore, we defined $N(G)$ to be $\gcd N(\rho)$ as $\rho$ varies over the representations of $G$.  We have proved:
\begin{equation}  \label{general}
\parbox{4.5in}{\emph{Let $G$ be split and almost simple. The trace $\Tr_\rho$ is zero for every representation $\rho$ of $G$  if and only if 
the characteristic of $F$ divides the integer $N(G)/E(G)$.}}
\end{equation}
We have now finished half of the proof of Theorem \ref{MT}\eqref{MT.nz}; it remains to check that the primes dividing $N(G)/E(G)$ are the primes in the right hand column of Table \ref{bad.primes}.
\end{borel*}

\section{The ratio $N(G)/E(G)$ for $G = \SL_n/\mmu{m}$}


In this section, we fix natural numbers $m$ and $n$ with $m$ dividing $n$, and we prove:
\begin{prop} \label{A.prop}
For $G = \SL_n/\mmu{m}$, the primes dividing $N(G)/E(G)$ are precisely the primes dividing
\[
\begin{cases}
\gcd(m,n/m) & \text{if $m$ is odd} \\
2\gcd(m,n/m) & \text{if $m$ is even.}
\end{cases}
\]
\end{prop}

Here $\mmu{m}$ denotes the group scheme of $m$-th roots of unity, identified with the corresponding scalar matrices in $\SL_n$.

In the important special cases where $G$ is simply connected ($m = 1$), $G$ is adjoint ($m = n$), or $n$ is square-free, the $\gcd$ in the proposition is 1, and we have that \emph{$N(G)/E(G)$ is $1$ if $m$ is odd and $2$ if $m$ is even}.  

\begin{lem} \label{A.E}
\[E(\SL_n / \mmu{m}) = \frac{m}{\gcd(m, n/m)}\,.\]
\end{lem}

\begin{proof}
Use the notation of \cite{BLie} for the simple roots and fundamental weights of the root system $A_{n-1}$ of $\SL_n$.  Let $\La$ denote the lattice generated by the root lattice $Q$ and
\[
\beta := \frac{n}{m}\,\omega_{n-1} = \frac{1}{m} \left( \alpha_1 + 2\alpha_2 + \cdots + (n-1)\alpha_{n-1} \right).
\]
We claim that $\La$ is identified with the cocharacter lattice $T_*$ for a pinning of $\SL_n/\mmu{m}$.  Certainly, $\La/Q$ is cyclic of order $m$, so it suffices to check that the set of inner products $\qform{\La, T^*}$ consists of integers.  But $T^*$ is the collection of weights $\sum c_i \omega_i$ with $c_i \in \Z$ such that $\sum_{i=1}^{n-1} ic_i$ is divisible by $m$. We have
$\qform{\beta, \sum c_i \omega_i} = \sum_i \frac{1}{m} ic_i$, which is an integer when $\sum c_i \omega_i$ is in $T^*$, so  $T_* = \La$ as claimed.

Finally, we compute:
\[
\qform{\beta, \alpha_{n-1}} = \frac{n}{m} \in \Z \quad \text{and} \quad \qform{\beta, \beta} = \qform{ \frac1m \sum i \alpha_i, \frac{n}{m} \omega_{n-1}} = \frac{n(n-1)}{m^2}.
\]
Since $m$ divides $n$, it is relatively prime to $n-1$, so the minimum multiplier of $\qform{\, , \,}$ that takes integer values on $T_*$ is $m/\!\gcd(m, n/m)$, as claimed.
\end{proof}

\begin{borel}{Weights of representations of $\SL_n / \mmu{m}$} \label{PGL.wts}
Fix the ``usual" pinning of $\SL_n$, where the torus $T$ consists of diagonal matrices and the dominant weights are the maps 
\[
\left( \begin{smallmatrix}
t_1& & \\ &\ddots & \\ && t_n \end{smallmatrix} \right) \mapsto \prod_{i=1}^{n-1} t_i^{e_i}
\]
where $e_1 \ge e_2 \ge \cdots \ge e_{n-1} \ge 0$.  Such a weight restricts to $x \mapsto x^{\sum e_i}$ on the center of $\SL_n$; in particular, $m$ divides $\sum e_i$ for every dominant weight $\la$ of a representation of $\SL_n/\mmu{m}$.  The proof of \cite[Lemma 11.4]{MG} shows that $m$ divides $N(W\la)$, hence \emph{$m$ divides $N(\SL_n/\mmu{m})$.}  
\end{borel}

\begin{borel*} \label{PGL.wts2}
We recall how to compute $N(W\la)$ from \cite[p.~136]{MG}.  Write $a_1 > a_2 > \cdots > a_{k-1} > a_k = 0$ for the distinct values of the exponents $e_i$ in $\la$, where $a_i$ appears $r_i$ times, so that $n = \sum r_i$.  
We have:
\begin{equation} \label{PGL.1}
N(W\la) = \frac{(n-2)!}{r_1!\, r_2! \cdots r_k!}  \left[ n\left( \sum_i r_i a_i^2 \right) - \left( \sum_i r_i a_i \right)^2 \right].
\end{equation}
The dominant weight $\la$ with $e_1 = m$ and $e_i = 0$ for $i > 1$ vanishes on $\mmu{m}$ and has $N(W\la) = m^2$ by \eqref{PGL.1}, so \emph{$N(\SL_n/\mmu{m})$ divides $m^2$}.
\end{borel*}

\begin{eg} \label{PGL.p2}
Let $\la$ be a dominant weight of $G$ and let $r_i, a_i$ be as in the preceding paragraph.  Suppose that 
\[
v_2\left(\sum r_i a_i \right) \ge v_2(n) > 0,
\]
where $v_2(x)$ is the 2-adic valuation of $x$, i.e., the exponent of the largest power of 2 dividing $x$.  We claim that
\begin{equation} \label{PGL.claim}
v_2(N(W\la)) > v_2(n).
\end{equation}

Write $\sum r_i a_i = 2^\theta t$ and $n = 2^\nu u$ where $\theta = v_2(\sum r_i a_i)$ and $\nu = v_2(n)$.  Our hypothesis is that $0 < \nu \le \theta$.  We rewrite \eqref{PGL.1} as:
\begin{equation} \label{PGL.2}
N(W\la) = \frac{(n-2)!}{r_1!\, r_2! \cdots r_k!} \left[ u \left( \sum_i r_i a_i^2 \right) - 2^{2\theta-\nu} t^2 \right]\cdot 2^\nu.
\end{equation}
Write $\ell$ for the minimum of $v_2(r_i)$, and fix an index $j$ such that $v_2(r_j) = \ell$.  Note that since $\sum r_i = n$, we have $\ell \le \nu \le 2\theta - \nu$.  

The first term on the right side of \eqref{PGL.2} has 2-adic value $\ge -\ell$ \cite[p.~137]{MG}.  The term in brackets has value $\ge \ell$.  Therefore, to prove claim \eqref{PGL.claim}, it suffices to consider the case where $v_2(\sum r_i a_i^2) = \ell$ and the first term on the right side of \eqref{PGL.2} has value $-\ell$; this latter condition implies that 
\begin{equation} \label{PGL.3}
s_2(n - 1) = s_2(r_1) + \cdots + s_2(r_{j-1}) + s_2(r_j - 1) + s_2(r_{j+1}) + \cdots + s_2(r_k),
\end{equation}
where $s_2$ denotes the number of 1's appearing in the binary representation of the integer \cite[p.~137]{MG}.  That is, when adding up the numbers 
$r_1, \ldots, r_{j-1}, r_j - 1, r_{j+1}, \ldots, r_k$ in base 2 (to get $n - 1$), there are no carries.  We check that this is impossible.

Suppose first that ${\ell}{\,<\,}{\nu}$.  Equation \eqref{PGL.3} implies that there are exactly two indices, say, $j, j'$ with $v_2(r_j) = v_2(r_{j'}) = \ell$.  As $2^{\ell + 1}$ divides $\sum r_i a_i$, it also divides $r_j a_j + r_{j'} a_{j'}$, hence $a_j$ and $a_{j'}$ have the same parity.  It follows that $2^{\ell + 1}$ divides $r_j a_j^2 + r_{j'} a_{j'}^2$, contradicting the hypothesis that $v_2(\sum r_i a_i^2) = \ell$.

We are left with the case where $\ell = \nu$.  
By \eqref{PGL.3}, $r_j$ is the unique $r_i$ with 2-adic valuation $\ell$.  As $v_2(\sum r_i a_i^2) = \ell$, the number $a_j$ is odd and we have:
$\ell = v_2\left(\sum r_i a_i\right) = \theta$.
Hence both $u\cdot(\sum r_i a_i^2)$ and $2^{2\theta-\nu}t$ have 2-adic valuation $\ell$.  It follows that the term in brackets in \eqref{PGL.2} has 2-adic valuation strictly greater than $\ell$, and claim \eqref{PGL.claim} is proved.
\end{eg}

\begin{proof}[Proof of Prop.~\ref{A.prop}]
We write $G$ for $\SL_n/\mmu{m}$.  Paragraphs \ref{PGL.wts} and \ref{PGL.wts2} give the bounds: $m$ divides $N(G)$ divides $m^2$.  Also, $N(G)$ divides $2n$ by Example \ref{Ad.eg}.  Applying Lemma \ref{A.E} gives:
\[
\text{$\gcd(m, n/m)$ divides $N(G)/E(G)$ divides $\gcd(m,n/m) \gcd(m, 2n/m)$.}
\]
This completes the proof for $m$ odd.

Clearly, an odd prime divides $N(G)/E(G)$ if and only if it divides $\gcd(m,n/m)$.  So suppose that $m$ is even and 2 does not divide $\gcd(m, n/m)$, i.e., $v_2(m) = v_2(n)$.  Then every dominant weight of a representation of $G$ satisfies the hypotheses of Example \ref{PGL.p2}, hence $v_2(N(G)) > v_2(n) = v_2(m)$.  By Lemma \ref{A.E}, $v_2(E(G)) = v_2(m)$, so 2 divides $N(G)/E(G)$.  This completes the proof of Prop.~\ref{A.prop}.
\end{proof}
\section{Conclusion of proof of Theorem \ref{MT}\eqref{MT.nz}} \label{conc.A2}

For a split and almost simple algebraic group $G$, we now verify that the primes dividing $N(G)/E(G)$ are those in the last column of Table \ref{bad.primes}.  Together with \eqref{general}, this will prove Theorem \ref{MT}\eqref{MT.nz}.

For $G$ simply connected, $E(G)$ is 1 and $N(G)$ is divisible precisely by the torsion primes of $G$, see \ref{dynk.def}.  We assume that $G$ is not simply connected and write $\Gt$ for the universal covering of $G$; obviously $N(\Gt)$ divides $N(G)$.

For $G = \PSp_{2n}$, $\SO_n$, or adjoint of type $E_7$, one combines Examples \ref{PSp.N} and \ref{E.C}; \ref{SO.N} and \ref{SO.E}; or \ref{E7.N} and \ref{E.sl}, respectively.


For $G$ adjoint of type $D_n$, we have $E(G) = 2$ by Example \ref{E.sl}.  Also, 4 divides $N(G)$ by \cite[15.2]{MG}.  On the other hand, the spinor representations of $\Gt$ have Dynkin index $2^{n-3}$ \cite{LaszloSorger}, and it is easy to use this as in Example \ref{E7.N} to construct a representation $\rho$ of $G$ with $N(\rho)$ a power of 2. This shows that $N(G)/E(G)$ is a power of 2 and is not 1.

Now let $G = \HSpin_{4n}$ for some $n \ge 3$.  The dual of the center of $\Spin_{4n}$ is the Klein four-group, and we write $\chi$ for the unique nonzero element that vanishes on the kernel of the map $\Spin_{4n} \ra \HSpin_{4n}$.  The gcd of $N(W\la)$ as $\la$ varies over the weights that restrict to $\chi$ (respectively, 0) on the center of $\Spin_{4n}$ is $2^{2n-3}$ (resp., divisible by 4) by \cite[p.~146]{MG}, hence $N(G)$ is a power of 2 and at least 4.  On the other hand, $E(\HSpin_{4n})$ is 1 or 2.  We conclude that $N(G)/E(G)$ is a power of 2 and is not 1.


For $G$ adjoint of type $E_6$, the number $N(G)$ is divisible by $N(\Gt) = 6$ and divides $2\hch = 24$ by Example \ref{Ad.eg}.    By Example \ref{E.sl}, $N(G)/E(G)$ is 2, 4, or 8.  This completes the proof of Theorem \ref{MT}\eqref{MT.nz}.$\hfill\qed$
%

\begin{eg} \label{SL.eg}
Suppose that the characteristic of $F$ is an odd prime $p$, and let $n$ be a natural number divisible by $p^2$. 
%
Every trace form of $\SL_n / \mmu{p}$ is zero by Theorem \ref{MT}\eqref{MT.nz}, even though the universal covering $\SL_n$ and adjoint group $\PGL_n$ have representations with nonzero trace forms. 
\end{eg}

\section{Proof of Theorem \ref{MT}\eqref{MT.nd}} \label{ndpf.sec}

We now prove Theorem \ref{MT}\eqref{MT.nd}; we show that the three statements
\begin{eqnarray} 
\parbox{4in}{\emph{The characteristic of $F$ is very good for $G$};} \label{H.vg} \\
\parbox{4in}{\emph{$\Lie(G)$ is a simple algebra and there is a representation $\rho$ of $G$ with $\Tr_\rho$ nonzero}; and} \label{H.nz} \\
\parbox{4in}{\emph{There is a representation $\rho$ of $G$ with $\Tr_\rho$ nondegenerate.}}\label{H.nd}
\end{eqnarray}
are equivalent.

Suppose \eqref{H.vg} holds.  Then $\Lie(G)$ is simple (as in \ref{Lie.simple}).  The existence of a representation $\rho$ with nonzero trace follows from Theorem \ref{MT}\eqref{MT.nz}, so \eqref{H.nz} holds.  
It is easy to check that for a representation $\psi$ of $\Lie(G)$, $\Tr_\psi([x,y],z) = \Tr_\psi(x,[y,z])$ for all $x, y, z \in \Lie(G)$.  So the radical of a trace form on $\Lie(G)$ is an ideal, and \eqref{H.nz} implies \eqref{H.nd}.

Now suppose that \eqref{H.vg} fails; we check that \eqref{H.nd} also fails.  By Theorem \ref{MT}\eqref{MT.nz}, we only need to consider those cases where the characteristic of $F$ appears in the middle column of Table \ref{bad.primes} and not in the right column, i.e., the cases:
\begin{enumerate}
\renewcommand{\theenumi}{\roman{enumi}}
\item \label{nd.sc} $G$ has type $G_2$ and $\chr F = 3$; or $G$ is $\Sp_{2n}$ and $\chr F = 2$; or  $G$ is $\SL_n / \mmu{m}$ and $\chr F$ is odd and divides $n/m$ but not $m$.
\item \label{nd.ad} $G$ is adjoint of type $E_6$ and $\chr F = 3$; or $G$ is $\SL_n / \mmu{m}$ and $\chr F$ divides $m$ but not $n/m$.
\end{enumerate}

We write $\pi \!: \Gt \ra G$ for the universal covering of $G$.
In case \eqref{nd.sc}, the kernel of $\pi$ is \'etale, so $\Lie(G)$ is a Weyl module by \ref{Lie.sc}.  For all three of the types listed, $\Lie(G)$ has a nontrivial submodule $M$, namely the subalgebra generated by the short roots (for $G_2$) or the center (in the other two cases).  It follows that $M$ is contained in the radical of $\Tr_\rho$---see e.g.\ \cite[6.2]{G:A1}---hence \eqref{H.nd} fails.

In case \eqref{nd.ad}, every representation $\rho$ of $G$ gives a representation $\rho\pi$ of $\Gt$ whose trace form $\Tr_{\rho\pi}$ vanishes on $\Lie(\Gt)$ by Theorem \ref{MT}\eqref{MT.nz} (for $E_6$) or \ref{PGL.wts} (for $\SL_n$).  Hence the image of $\dpi$ is a totally isotropic subspace for $\Tr_\rho$.  As
\[
\dim (\im \dpi) = \dim \Gt - \dim(\ker \dpi) = \dim G - 1
\]
is strictly greater than half the dimension of $G$, the form $\Tr_\rho$ is degenerate and \eqref{H.nd} fails.  This concludes the proof of Theorem \ref{MT}\eqref{MT.nd}.$\hfill\qed$

\section{Richardson's condition}

In the literature, the weak version of the ``if" direction of Theorem \ref{MT}\eqref{MT.nd} is used to 
deduce ``Richardson's condition" from \cite[p.~3]{Richardson}.  Our slightly finer version of the ``if" direction gives a slightly finer version of Richardson's condition; we state it here for the convenience of the reader.   As in Theorem \ref{MT}, $G$ is a split almost simple algebraic group over a field $F$.

\begin{prop}
If the characteristic of $F$ is very good for $G$, then there is a representation $\rho \!: G \ra \GL(V)$ such that $\drho$ is an injection $\Lie(G) \injects \gl(V)$ and there is a subspace $M$ of $V$ such that $V = \drho(\Lie(G)) \oplus M$, $\Id_V$ is in $M$, and $\Ad(\rho(G)) M \subseteq M$.
\end{prop}

\begin{proof}
Theorem \ref{MT}\eqref{MT.nd} gives a representation $\rho$ so that $\Tr_\rho$ is nondegenerate.  In particular, the restriction of the symmetic bilinear form $(x, y) \mapsto \trace(xy)$ on $\gl(V)$ to $\drho(\Lie(G))$ is nondegenerate.  (And obviously $\drho$ must be injective.)  

Take $M$ to be the space of $x \in \gl(V)$ such that $\trace(\drho(\Lie(G))x) = 0$.  Trivially, $M$ is invariant under $\Ad(\rho(G))$.  Nondegeneracy of $\Tr_\rho$ shows that $M$ meets $\drho(\Lie(G))$ only at 0, and dimension count shows that $V = \drho(\Lie(G)) \oplus M$.  As $G$ is semisimple, the image $\rho(G)$ is contained in $\SL(V)$, hence $\drho(\Lie(G))$ lies in $\sl(V)$; i.e., $\Id_V$ belongs to $M$.
\end{proof}

The proposition is essentially known, but the usual argument as in \cite[\S5]{Richardson},  \cite[2.6]{Jantzen:nil}, \cite[p.~48]{Hum:cc}, or \cite[p.~184]{SpSt} is different.  For example, the usual approach to treating an adjoint group $G$ of type $C_n$ or $D_n$ replaces $G$ with its covering $G' = \Sp_{2n}$ or $\SO_{2n}$ and then gives a representation of $G'$ with the desired properties.

\section{Complements: characteristic 2}

In characteristic 2, one might prefer to consider the quadratic form
\[
s_\rho \!: x \mapsto -\trace \left( \wedge^2 \drho(x) \right)
\]
instead of the symmetric bilinear form $\Tr_\rho$.  The form $s_\rho$ 
gives the negative of the ``degree 2" coefficient of the characteristic polynomial of $\drho(x)$.  (Because $\drho(\g)$ consists of trace zero matrices, $s_\rho$ is the map $x \mapsto \trace(\drho(x)^2)/2$; our definition has the advantage that it obviously makes sense also in characteristic 2.)  The bilinear form derived from $s_\rho$---i.e., $(x,y) \mapsto s_\rho(x+y) - s_\rho(x) - s_\rho(y)$---is $\Tr_\rho$.

Theorem \ref{MT}\eqref{MT.nz} is easy to extend.  In case $G$ is simply connected, $\Lie(G)$ is a Weyl module by \ref{Lie.sc} and $s_\rho$ is zero if and only if $\Tr_\rho$ is zero by \cite[Prop.~6.4(1)]{G:A1}.  That is, the conditions in Theorem \ref{MT}\eqref{MT.nz} are equivalent to: \emph{For every representation $\rho$ of $G$, the quadratic form $s_\rho$ is zero.}  (This is true in all characteristics but is only nontrivial in characteristic 2.)

Alternatively, one can proceed as follows.  The bilinear form $\bt$ on $\Lie(\Gt_\Z)$ is even \cite[Prop.~4]{GrossNebe}, so it is the bilinear form derived from a unique quadratic form $\qt$ on $\Lie(\Gt_\Z)$.  The form $\qt$ extends to a rational-valued quadratic form on $\Lie(G_\Z)$ and we write $E_q(G)$ for the smallest positive rational number such that $E_q(G)\,\qt$ is integer-valued on $\Lie(G_\Z)$.  The number $E_q(G)$ is $E(G)$ or $2E(G)$, and both cases can occur.  (E.g., take $G = \Gt$ or $\SO_{2\ell}$, respectively.)  The statements and proofs of \eqref{single} and \eqref{general} go through if we replace $\Tr_\rho$, $E(G)$, and $b$ with $s_\rho$, $E_q(G)$, and $E_q(G)\,\qt$ respectively.

\section{Complements: non-split groups} \label{nonsplit}

We can extend our results above to the case where $G$ is not split, i.e., we can replace the hypotheses ``$G$ is split and almost simple" with ``$G$ is absolutely almost simple".  Indeed, suppose that $G$ is absolutely almost simple over $F$, i.e., there is a split  and almost simple group $G'$ over $F$ and an isomorphism $f \!: G' \ra G$ defined over a separable closure $\Fsep$ of $F$. 
Fix a pinning for $G'$ and write $b'$ for the indivisible bilinear form on $\Lie(G'_\Z)$ defined in \ref{b.def}.  Clearly, the automorphism group of $G'$---which is an affine group scheme over $\Z$---leaves $b'$ invariant, so it maps into the orthogonal group of $b'$.  Galois descent (via $f$) gives a $G$-invariant symmetric bilinear form $b_F$ on $\Lie(G)$ such that the differential $\df$ identifies $b'_F \ot \Fsep$ with $b_F \ot \Fsep$.

Given a representation $\rho$ of $G$ over $F$, we get a representation $\rho f$ of $G'$ over $\Fsep$ and an integer $N(\rho f)$ defined in \ref{Nrho.def}; put $N(\rho) := N(\rho f)$.  (In the special case where $G$ is split over $F$, this agrees with our previous definition.)  We define $N(G)$ as in \eqref{NG.def}; it is the gcd of $N(\rho)$ as $\rho$ varies over the representations of $G$ defined over $F$.  Obviously, $N(G)$ is divisible by $N(G_{\Fsep})$---i.e., $N(G')$---and it depends on the field $F$.  

We put $E(G) := E(G')$.  It does not depend on the field $F$.

With these definitions for $N(G)$ and $E(G)$, conclusions (1) and (2) of Proposition \ref{divides} hold for absolutely almost simple $G$.  Indeed, it suffices to check them over $\Fsep$, where they hold by the original version of the proposition.  It follows immediately that the conclusions of \eqref{single} and \eqref{general} hold for every absolutely almost simple group $G$.  

We now extend Theorem \ref{MT}.  Recall that there is a natural action of the absolute Galois group $\Gal(F)$ of $F$ on the Dynkin diagram $\D$ of $G$ \cite[2.3]{Ti:Cl}.  As in \cite[p.~54]{Ti:Cl}, we say, for example, that $G$ \emph{has type $^3\!D_4$} if $\D$ has type $D_4$ and the image of the map $\Gal(F) \ra \Aut(\D)$ has order 3.    We say that the characteristic of $F$ is not very good for $G$ if and only if it is not very good for the corresponding split group $G'$; these primes are listed in the middle column of Table \ref{bad.primes}.

\begin{specthmp} \label{MTp}
Let $G$ be an absolutely almost simple algebraic group over a field $F$.
\begin{enumerate}
\item Every representation $\rho$ of $G$ over $F$ has $\Tr_\rho$ degenerate if and only if the characteristic of $F$
\[
\begin{cases}
\text{divides $2n$} & \text{if $G$ has type $^2\!A_{n-1}$ for some odd $n \ge 3$;} \\
\text{is $2$ or $3$} &\text{if $G$ has type $^3\!D_4$ or $^6\!D_4$;} \\
\text{is not very good for $G$} & \text{otherwise.}
\end{cases}
\]
\item Suppose $G$ is not simply connected and not of type $A$.  Every representation of $G$ has $\Tr_\rho$ zero if and only if the characteristic of $F$ is as in the table:
\[
\begin{tabular}{c|c}
type of $G$&$\chr F$ \\ \hline
\upstrut{3} $B_n$ $(n \ge 3)$; $C_n$ $(n \ge 2)$; $^1\!D_n$ or $^2\!D_n$ $(n \ge 4)$; or $E_6$ & $2$ \\
$^3\!D_4$, $^6\!D_4$, $E_7$&$2$ or $3$
\end{tabular}
\]
\end{enumerate}
\end{specthmp}

Regarding the omitted cases in part \eqref{MT.nz}, for $G$ simply connected, the number $E(G)$ is 1, so every representation $\rho$ of $G$ has $\Tr_\rho$ zero if and only if the characteristic divides $N(G)$ by \eqref{general}; this number (using that $G$ is simply connected) is calculated in \cite[\S\S11--16]{MG}.  We leave the type $A$ case of \eqref{MT.nz} as an exercise for the reader.

\begin{proof}[Proof of Theorem \ref{MTp}]
To prove \eqref{MT.nz}, by \eqref{general} it remains to show that the primes in the table are those dividing $N(G)/E(G)$.  As $N(G')$ divides $N(G)$ and $E(G')$ equals $E(G)$, we have the trivial equation 
\begin{equation} \label{MTp.1}
\frac{N(G)}{E(G)} = \frac{N(G)}{N(G')} \frac{N(G')}{E(G')}
\end{equation}
where all three terms are integers.  The primes dividing $N(G')/E(G')$ are listed in Table \ref{bad.primes}, so it suffices to check which primes divide $N(G)/N(G')$ and are not in that table.

For $G$ adjoint of type $E_6$, the proof that $N(G)/E(G)$ is a power of 2 from the end of \S\ref{conc.A2} goes through without change.

The proof of \cite[10.11]{MG} shows that every prime dividing $N(G)/N(G')$ divides the exponent of $P/Q$ (the weight lattice modulo the root lattice) or the order of the image of $\Gal(F) \ra \aut(\D)$.  For $G$ of type $B_n$ ($n \ge 3$); $C_n$ ($n \ge 2$); $^1\!D_n$ or $^2\!D_n$ ($n \ge 4$); or $E_7$, the exponent of $P/Q$ is 2 and the image of $\Gal(F) \ra \aut(\D)$ has order at most 2.  As 2 divides $N(G')/E(G')$, part \eqref{MT.nz} is proved for these groups.

For $G$ adjoint of type $^3\!D_4$ or $^6\!D_4$, write $\Gt \ra G$ for the universal covering of $G$.  The number $N(\Gt)$ is 6 or 12 by \cite[16.5]{MG} and divides $N(G)$.  As $E(G)$ is 2 by Example \ref{E.sl}, $N(G)/E(G)$ is divisible by 3.  Part \eqref{MT.nz} of the theorem is proved.

(We remark that applying the argument from the two previous paragraphs in the case where $G$ has type $A_{n-1}$ shows that every prime dividing $N(G)/N(G')$ divides $2n$.  If $n$ is odd and $\ge 3$ and $G$ has type $^2\!A_{n-1}$, then 2 divides $N(\Gt)$ by \cite[12.6]{MG} hence also $N(G)$, yet $E(G)$ is odd by Lemma \ref{A.E}, so $N(G)/E(G)$ is even.)


\newcommand{\Hvgp}{(\ref{H.vg}$^\prime$)}
\smallskip
We now prove part \eqref{MT.nd} by imitating \S\ref{ndpf.sec}.  We replace \eqref{H.vg} with the condition that the characteristic of $F$ is not as in the statement of Theorem \ref{MTp}\eqref{MT.nd}; we denote this condition by \Hvgp. 

Suppose that \Hvgp\ holds.  The characteristic is very good for $G$ and $\Lie(G) \ot \Fsep$ is simple as in \ref{Lie.simple}, hence $\Lie(G)$ is simple.  If $G$ is neither simply connected nor of type $A$, then there is a representation $\rho$ of $G$ with $\Tr_\rho$ nonzero by part \eqref{MT.nz} and \eqref{H.nz} holds.  If $G$ is simply connected, then checking \cite{MG} verifies that $N(G)$ is not divisible by the characteristic and again \eqref{H.nz} holds.  In the remaining case where $G$ has type $A$, the characteristic does not divide $N(G')$ by \ref{PGL.wts2} nor does it divide $N(G)/N(G')$ by the discussion above; by \eqref{MTp.1}, we find that \eqref{H.nz} holds.

As in \S\ref{ndpf.sec}, \eqref{H.nz} trivially implies \eqref{H.nd}.

Finally, suppose that \Hvgp\ fails; we will show that \eqref{H.nd} fails.  We assume that the characteristic is very good, otherwise \eqref{H.nd} fails because it does so over $\Fsep$.  That is, we are in one of the cases
\begin{enumerate}
\renewcommand{\theenumi}{\roman{enumi}}
\item $\chr F = 3$ and $G$ has type $^3\!D_4$ or $^6\!D_4$; or 
\item $\chr F = 2$ and $G$ has type $^2\!A_{n-1}$ for some odd $n$.
\end{enumerate}
But in these cases the characteristic divides $N(G)/E(G)$ by the proof of part \eqref{MT.nz} above, and \eqref{H.nd} fails.
\end{proof}

%

\begin{eg}
Let $F$ be a field of prime characteristic $p$ with a central division $F$-algebra $A$ of degree $p$.  Take $G$ to be the group $\SL(A)$ whose $F$-points are the elements of $A$ with determinant 1.  This group is simply connected, so $N(G)/E(G)$ is $p$ by \cite[11.5]{MG}.  That is, $\Tr_\rho$ is zero for every representation $\rho$ of $G$ over $F$.  On the other hand, $N(G_{\Fsep})$ is 1, so there are representations of $G$ defined over $\Fsep$ (e.g., the natural representation of $\SL_n$) that have a trace form that is not zero.

A similar statement holds for groups of type $^3\!D_4$ or ${^6\!D_4}$ over fields of characteristic 3.
\end{eg}


\section{Trace forms and Lie algebras} \label{trace.sec}

This section collects some results regarding $\Tr_\psi$, where $\psi$ is a representation of the Lie algebra of an algebraic group $G$ and we do not assume that $\psi$ is the differential of a representation of $G$.

\begin{borel*} \label{block.hyp}
Fix a positive integer $n$ and assume that the characteristic of $F$ is a prime dividing $n$ and $\ne 2, 3$.
The Lie algebra $\sl_n$ of trace zero $n$-by-$n$ matrices has center $\c$ the scalar matrices and $\sl_n / \c$ is simple \cite[2.6]{St:aut}.  We give a new proof of:
\end{borel*}

\begin{prop}[Block \protect{\cite[Th.~6.2]{Block:trace}}]  \label{Block.cor}
Under the hypotheses of \ref{block.hyp}, 
every representation of $\sl_n / \c$ has zero trace form.
\end{prop}

\begin{proof}
For sake of contradiction, suppose that there is an irreducible representation $\psi$ of $\sl_n/\c$ with nonzero trace form.
Then $\psi$ is restricted by \cite[Th.~5.1]{Block:trace} (using that $F$ has characteristic $\ne 2, 3$).  The composition of $\psi$ with $\sl_n \ra \sl_n/\c$ is  a restricted irreducible representation of $\sl_n$, which is the differential of a representation $\rho$ of $\SL_n$ by \cite{Curtis} and \cite{St:rep}.  

By construction $\Tr_\rho$ is not zero and $\drho$ vanishes on the scalar matrices.  
Identifying the center of $\SL_n$ with the (non-reduced) group scheme $\mmu{n}$ identifies the restriction of $\rho$ to $\mmu{n}$ with a map $x \mapsto x^\ell$.  Our hypothesis on $\drho$ says that $\ell$ is divisible by the characteristic $p$ of $F$, hence $\rho$ factors through the natural map $\SL_n \ra \SL_n/\mmu{p}$.  Paragraph \ref{PGL.wts} says that $N(\rho)$ is divisible by $p$, hence $\Tr_\rho$ vanishes by \eqref{single}, a contradiction.

Since every irreducible representation has zero trace form, the same holds for every representation like at the end of the proof of Proposition \ref{divides}.
\end{proof}

\begin{boreli}{Proof of Theorem \ref{alg}} \label{alg.pf}
Let $G$ be a group of type $E_8$, and suppose that there is a representation $\psi$ of $\Lie(G)$ such that $\Tr_\psi$ is not zero.  We may assume that $\psi$ is irreducible.  Then Theorem \ref{Premet} implies that $\psi$ is restricted, hence is the differential of a representation of $G$.  Theorem \ref{MT}\eqref{MT.nz} implies that the characteristic of $F$ is $\ne 2, 3, 5$.
$\hfill\qed$
\end{boreli}

We close by proving that over a field of characteristic 5, the Lie algebra of a group of type $E_8$ ``has no quotient trace form".  For a Lie algebra $L$ over $F$ and a representation $\psi$ of $L$, write $\rad \psi$ for the radical of the trace bilinear form $\Tr_\psi$; it is an ideal of $L$.  We prove:
\begin{cor} \label{E8.cor}
For every representation $\psi$ of every Lie algebra $L$ over a field of characteristic $5$, the quotient $L/\rad \psi$ is \emph{not} isomorphic to the Lie algebra of an algebraic group of type $E_8$.
\end{cor}

\begin{proof}
Suppose the corollary is false.  That is, suppose that there is a group $G$ of type $E_8$ and a Lie algebra $L$ with a representation $\psi$ and a surjection $\pi \!: L \ra \g$ with kernel the radical of $\Tr_\psi$.

By \cite[Lemma 2.1]{Block:trace}---using that the characteristic is $\ne 2, 3$---we may assume that the radical of $\Tr_\psi$ is contained in the center of $L$, i.e., $L$ is a central extension of $\g$.  It follows that there is a map $f \!: \g \ra L$ such that $\pi f$ is the identity \cite[Th.~6.1(c)]{St:gen}.  Clearly, the representation $\psi f$ of $\g$ has nonzero trace form.  
As in the proof of Proposition \ref{Block.cor}, we deduce that $G$ has a representation $\rho$ such that $\Tr_\rho$ is not zero, but this is impossible by Theorem \ref{MT}\eqref{MT.nz}.
\end{proof}

\section{Appendix: On trace forms of Lie algebras of type $E_8$\\{\protect\textit{By Alexander Premet}}}


\newcommand{\te}{\mathfrak{t}}
\newcommand{\es}{\mathfrak{s}}
\newcommand{\el}{\mathfrak{l}}
\newcommand{\z}{\mathfrak{z}}

\renewcommand{\g}{\mathfrak{g}}

All basic notions and results of modular Lie theory used in this
appendix can be found in \cite{Premet:KW} and references therein.

Let $G$ be an algebraic group of type $E_8$ over an algebraically
closed field of characteristic $p>0$ and $\g={\rm Lie}(G)$. It is
well known that $\g$ is a simple Lie algebra carrying an
$(\Ad G)$-equivariant $[p]$-th power map $x\mapsto x^{[p]}$. Since
the universal enveloping algebra $U(\g)$ is a finite module over its
central subalgebra generated by all $x^p-x^{[p]}$ with $x\in\g$, all
irreducible $\g$-modules are finite dimensional. Furthermore, for
every irreducible $\g$-module $M$ there is a linear function
$\xi=\xi_M$ on $\g$ such that $x^p-x^{[p]}$ acts on $M$ as the
scalar operator $\xi(x)^p\Id_M$. The function 
$\xi_M$ is
called the $p$-{\it character} of $M$. Denote by $I_\xi$ the
two-sided ideal of $U(\g)$ generated by all elements
$x^p-x^{[p]}-\xi(x)^p$, where $x\in\g$. The factor-algebra of
$U(\g)/I_\xi$ is called the {\it reduced enveloping algebra}
associated with $\xi$ and denoted $U_\xi(\g)$. It has dimension
$p^{\dim\,\g}$. Clearly, $M$ is a $U_\xi(\g)$-module. We say that
$M$ is {\it restricted} if $\xi_M=0$.

For $p>3$, the theorem below was first proved by Richard Block in
\cite{Block:trace}. The aim of this appendix is to give a proof valid in any
positive characteristic.
\begin{thm} \label{Premet}
If $\psi \!: \g\rightarrow\,\mathfrak{gl}(V)$ is an irreducible
representation with $\Tr_\psi\ne 0$, then $V$ is a
restricted $\g$-module.
\end{thm}

\begin{proof}
Suppose $\psi$ is not restricted and let $\chi$ be the $p$-character
of $V$. Then $\chi$ is a nonzero linear function on $\g$.  We show that $\Tr_\psi$ is zero.

Let $T$ be a maximal torus of $G$ and $\te:={\rm Lie}(T)$. As in sections \ref{N.sec} and \ref{Lie.sec},
we write $R$ for the root system of $G$ relative
to $T$ and $h_\alpha$ for the image of the coroot
$\alpha^\vee$ in $\te:={\rm Lie}(T)$. (In our case, the group $G$ is
both adjoint and simply connected.) Since $\g$ is an irreducible
$(\Ad G)$-module, every nonzero adjoint $G$-orbit spans $\g$.
Thus, replacing $\te$ by its $G$-conjugate if necessary, we may
assume that $\chi(h_\beta)\ne 0$ for some $\beta\in R$.

There are root vectors $e_{\pm\beta}\in\g^{\pm\beta}$ such that
$\es:=Fe_{-\beta}\oplus Fh_\beta\oplus Fe_\beta$ is isomorphic to
$\mathfrak{sl}_2$. Replacing $\te$ by its conjugate
$(\Ad x_{-\beta}(\lambda))(\te)$ for a suitable
$x_{-\beta}(\lambda)$ in the unipotent root subgroup $U_{-\beta}$ of $G$,
we may assume without loss of generality that $\chi\vert_\es\ne
0$ and $\chi(e_\beta)=0$.
Then every $\es$-composition factor $M$ of $V$ is a baby Verma
module, that is, $M\cong Z_\xi(a)$, where $\xi=\chi|_\es$ and $a\in
F$ is a root of the equation $X^p-X=\xi(h_\alpha)^p$. Note that $\dim
M=p$, the operator $h_\beta$ acts semisimply on $M$, and the
$h_\beta$-weights of $M$ are $a$, $a-2, \ldots,  a-2(p-1)$.

First suppose $p>3$. Then ${\rm
trace}_M(h_\beta^2)=\textstyle{\sum}_{i=0}^{p-1}\,(a-2i)^2=pa^2-2ap(p-1)+\frac{2}{3}(p-1)p(2p-1)=0.$
Since this holds for every $\es$-composition factor $M$ of $V$, we
obtain $\Tr_\psi(h_\beta,h_\beta)=0$.  As $\g$ is a simple Lie algebra and $\Tr_\psi$ is $\g$-invariant, $\Tr_\psi$ is a multiple of the form $b_F$ from section \ref{crux.sec}.  Hence $\Tr_\psi$ is zero.
%

Next suppose $p=3$. Then the $h_\beta$-weights of $M$ are $a$, $a+1$,
$a-1$, hence $\trace_M(h_\beta^2)=a^2+(a+1)^2+(a-1)^2=2$. It
follows that $\trace_M(h_\beta^2)$ is independent of $M$. Since
all $\es$-composition factors of $V$ are $3$-dimensional, we deduce
that $\Tr_\psi(h_\beta,h_\beta)=2(\dim V)/3$. Note that $\es$
can be included into a Levi subalgebra of type $A_7$; call it $\el$.
Since $\es \subset \el$, all $\el$-composition factors of $V$ have
the same nonzero $p$-character. But then the Kac--Weisfeiler
conjecture (which holds for $\mathfrak{sl}_8$ in characteristic $3$
thanks to \cite[Th.~3.10]{Premet:KW}) implies that all such factors have
dimension divisible by $9$. Then $9$ divides $\dim\, V$, forcing ${\rm
Tr}_\psi(h_\beta,h_\beta)=0$. As in the $p > 3$ case, $\Tr_\psi$ is zero.

Finally, suppose $p=2$. Then the $\mathfrak{sl}_2$-algebra
$\es=Fe\oplus Fh\oplus Ff$ is nilpotent and $h$ lies in the center
of $\es$. However, the reduced enveloping algebra $U_\xi(\es)$ is
semisimple whenever $\xi(h)\ne 0$. Indeed, $U_\xi(\es)$ then
possesses two non-equivalent $2$-dimensional irreducible modules,
$M$ and $N$, induced from $1$-dimensional modules over a Borel
subalgebra of $\es$. The central element $h$ of $\es$ acts on $M$
and $N$ by different scalars. There are exactly two choices here,
namely, $a$ and $a+1$, where $a$ is a root of the equation
$X^2-X=\xi(h)^2$. As a consequence, $U_\xi(\es)$ maps onto a direct
sum of two copies of ${\rm Mat}_2(F)$. Since $\dim
U_\xi(\es)=2^3=8$, this map is an isomorphism. Thus, $U_\xi(\es)$ is
semisimple with two isoclasses of simple modules, both of which are
$2$-dimensional.

Suppose now that we have found two commuting
$\mathfrak{sl}_2$-subalgebras $\es_i=Fe_i\oplus Fh_i\oplus Ff_i$ in
$\g$, where $i=1,2$, such that
\begin{enumerate}
\renewcommand{\theenumi}{\alph{enumi}}
\item the sum $\es_1+\es_2$ is direct;

\item $\chi(h_i)\ne 0$ for $i=1,2$;

\item  $e_1\in\g^\gamma$ and $f_1\in\g^{-\gamma}$ for some $\gamma\in R$.
\end{enumerate}
Our preceding remark then would show that $V$ is a semisimple module
over the subalgebra $U_\chi(\es_1\oplus \es_2)\cong
U_\chi(\es_1)\otimes U_\chi(\es_2)$ of $U_\chi(\g)$ (to ease
notation, we identify $\chi$ with its restriction to $\es_i$,
$i=1,2$). Let $M$ and $N$ be two irreducibles for $U_\chi(\es_1)$
described earlier. Then $V$ decomposes as a tensor product
$V=(M\otimes P)\bigoplus (N\otimes Q)$ for some semisimple
$U_\chi(\es_2)$-modules $P$ and $Q$. Therefore,
\[
\Tr_\psi(e_1,f_1)=r\dim P +s\dim Q,
\]
where $r=\trace_{M}(e_1f_1)$ and $s=\trace_N(e_1f_1)$. As
both $P$ and $Q$ must have even dimension by our preceding remark,
this would yield $\Tr_\psi(\g^\gamma,\g^{-\gamma})=0$.  Hence $\Tr_\psi = 0$ by equation \eqref{xlong}, using that $\gamma$ is (trivially) a long root.

So it remains to find two commuting $\mathfrak{sl}_2$-triples as
above. We adopt Bourbaki's numbering of simple roots; see \cite{BLie}.
Since $\chi\ne 0$ and the adjoint $G$-orbit of $e_{\alpha_7}$ spans
$\g$ by the simplicity of $\g$, we may assume that
$\chi(e_{\alpha_7}) \ne 0$. If $\chi(h_{\alpha_6})\ne 0$ and
$\chi(h_{\alpha_8})\ne 0$, then we can take
$\es_1=Fe_{\alpha_6}\oplus Fh_{\alpha_6}\oplus Fe_{-\alpha_6}$ and
$\es_2=Fe_{\alpha_8}\oplus F h_{\alpha_8}\oplus Fe_{-\alpha_8}$. If
this is not the case, then we replace $\te$ by $(\Ad x_{\alpha_7}(t_0))(\te)$ for a suitable $x_{\alpha_7}(t_0)$ in
the unipotent root subgroup $U_{\alpha_7}$ of $G$.

There exists $b\in F^\times$ such that for every $t\in F$ we have
\[
(\Ad x_{\alpha_7}(t))(h_{\alpha_i})=h_{\alpha_i}+tb[e_{\alpha_7},
h_{\alpha_i}]=h_{\alpha_i}+tbe_{\alpha_7},\quad i=6,8.
\]
As
$\chi(e_{\alpha_7})\ne 0$, we can find $t_0\in F$ such that 
$\chi(h_{\alpha_6}+t_0be_{\alpha_7})\ne 0$ and
$\chi(h_{\alpha_8}+t_0be_{\alpha_7})\ne 0$. Hence we can take
$(\Ad x_{\alpha_7}(t_0))(\es_i)$, $i=1,2$, as our
$\mathfrak{sl}_2$-triples. This completes the proof.
\end{proof}


\providecommand{\bysame}{\leavevmode\hbox to3em{\hrulefill}\thinspace}
\providecommand{\MR}{\relax\ifhmode\unskip\space\fi MR }
\providecommand{\MRhref}[2]{%
  \href{http://www.ams.org/mathscinet-getitem?mr=#1}{#2}
}
\providecommand{\href}[2]{#2}

\end{document}

%% file: skip-thm.tex
\usepackage{amsthm}

\newtheorem{thm}[equation]{Theorem}
\newtheorem{lem}[equation]{Lemma}
\newtheorem{cor}[equation]{Corollary}
\newtheorem{prop}[equation]{Proposition}

\newtheorem*{thm*}{Theorem}
\newtheorem*{prop*}{Proposition}
\newtheorem*{cor*}{Corollary}
\newtheorem*{lem*}{Lemma}
\newtheorem*{MT*}{Main Theorem}

\newtheorem{specthm}{Theorem}

\theoremstyle{definition} %
\newtheorem{defn}[equation]{Definition}
\newtheorem*{defn*}{Definition}

\newtheorem{eg}[equation]{Example}

\theoremstyle{remark} %

\newtheorem*{rmk*}{Remark}
\newtheorem*{rmks*}{Remarks}

\newtheoremstyle{exercise}
  {3pt}
  {3pt}
  {\small}
  {}
  {\sc\small}
  {.}
  {.5em}
   {}     
  {}

\theoremstyle{exercise}

\makeatletter
{\renewcommand{\theequation}{#1}}%
{\renewcommand{\theequation}{\arabic{equation}}\addtocounter{equation}{-1}\global\@ignoretrue}
\makeatother

\makeatletter
{\renewcommand{\theequation}{#1}\begin{eqnarray}}%
{\end{eqnarray}\renewcommand{\theequation}{\arabic{equation}}\addtocounter{equation}{-1}\global\@ignoretrue}
\makeatother

\makeatletter
\newenvironment{borel}[1]%
{\smallskip \refstepcounter{equation}\noindent{\textbf{\theequation.} }{{\textbf{#1.}}}}%
{\smallskip \global\@ignoretrue}
\makeatother

\makeatletter
{\smallskip \refstepcounter{equation}\noindent{\textbf{\theequation.} }{{\textbf{#1}}}}%
{\smallskip \global\@ignoretrue}
\makeatother

\makeatletter
{\smallskip \refstepcounter{equation}{\sc \theequation}{\sc (#1).}}%
{\smallskip \global\@ignoretrue}
\makeatother

\makeatletter
\newenvironment{boreli}[1]%
{\smallskip\refstepcounter{equation}\noindent{\textbf{\theequation.}}{\textsl{ #1.}}}%
{\smallskip\global\@ignoretrue}
\makeatother

\makeatletter
\newenvironment{borel*}%
{\smallskip \refstepcounter{equation}\noindent{\textbf{\theequation.}}}%
{\global\@ignoretrue}
\makeatother

\newcommand{\flist}[1]{\hangindent\leftmargini\textup{(1)}\hskip\labelsep {#1}%
\begin{enumerate}%
\setcounter{enumi}{1}%
}

%% file: qformsl-def.tex
\newcommand{\ot}{\otimes}

\newcommand{\eand}{\quad\text{and}\quad}

\DeclareRobustCommand{\upstrut}[1]{\rule{0mm}{#1 ex}}

\DeclareMathOperator{\Lie}{Lie}

\newcommand{\C}{{\mathbb{C}}}        
\newcommand{\Q}{{\mathbb{Q}}}        
\newcommand{\Z}{{\mathbb{Z}}}        

\newcommand{\Zm}[1]{\Z/{#1}\Z}

\newcommand{\La}{\Lambda}
\newcommand{\la}{\lambda}
\newcommand{\D}{\Delta}

\newcommand{\oddots}{{\mathinner{\mkern1mu\raise1pt\vbox{\kern7pt\hbox{.}}\mkern2mu\raise4pt\hbox{.}\mkern2mu\raise7pt\hbox{.}\mkern1mu}}}

\newcommand{\Fsep}{F_{{\mathrm{sep}}}}

\newcommand{\qform}[1]{{\left\langle{#1}\right\rangle}}                   

\newcommand{\basemu}{\boldsymbol{\mu}}
\newcommand{\mmu}[1]{\basemu_{#1}}     

\DeclareMathOperator{\Spin}{Spin}           
\newcommand{\Sp}{\mathrm{Sp}}
\DeclareMathOperator{\SL}{SL}
\DeclareMathOperator{\GL}{GL}
\DeclareMathOperator{\PGL}{PGL}
\newcommand{\SO}{\mathrm{SO}}

\DeclareMathOperator{\Ad}{Ad}

\DeclareMathOperator{\Tr}{Tr}

\DeclareMathOperator{\Gal}{Gal}  

\DeclareMathOperator{\im}{im}

\DeclareMathOperator{\chr}{char}

\DeclareMathOperator{\Id}{Id}

\DeclareMathOperator{\aut}{Aut}
\DeclareMathOperator{\Aut}{Aut}

\newcommand{\Hom}{{\mathrm{Hom}}}

\newcommand{\iso}{\xrightarrow{\sim}}
\newcommand{\injects}{\hookrightarrow}

\newcommand{\ra}{\rightarrow}

%% file: vanish.bbl
\begin{thebibliography}{MPR90}

\bibitem[Bl]{Block:trace}
R.E. Block, \emph{Trace forms on {L}ie algebras}, Canad. J. Math. \textbf{14}
  (1962), 553--564.

\bibitem[BlZ]{BlockZ}
R.E. Block and H.~Zassenhaus, \emph{The {L}ie algebras with a nondegenerate
  trace form}, Illinois J. Math. \textbf{8} (1964), 543--549.


\bibitem[Bo]{BLie}
N.~Bourbaki, \emph{Lie groups and {L}ie algebras: Chapters 4--6},
  Springer-Verlag, Berlin, 2002.
  
\bibitem[BR]{BR}
P. Bardsley and R.W. Richardson, \emph{\'Etale slices for algebraic transformation groups in characteristic $p$}, Proc. London Math. Soc. (3) \textbf{51} (1985), 295--317.  
  
\bibitem[Ca]{Carter:big}
R.W. Carter, \emph{Finite groups of Lie type: conjugacy classes and complex characters}, Wiley-Interscience, 1985.
  
\bibitem[Cu]{Curtis}
C.W. Curtis, \emph{Representations of {L}ie algebras of classical type with
  applications to linear groups}, J. Math. Mech. \textbf{9} (1960), 307--326.
  
\bibitem[DG]{DG}
M. Demazure and P. Gabriel, \emph{Groupes alg\'ebriques}, North-Holland, Amsterdam, 1970.  

\bibitem[Dy]{Dynk:ssub}
E.B. Dynkin, \emph{Semisimple subalgebras of semisimple {L}ie algebras}, Amer.
  Math. Soc. Transl. (2) \textbf{6} (1957), 111--244, [Russian original: Mat.\
  Sbornik N.S.\ \textbf{30(72)} (1952), 349--462].

\bibitem[Ga]{G:A1}
S.~Garibaldi, \emph{Orthogonal representations of twisted forms of ${\SL}_2$}, Representation Theory \textbf{12} (2008), 453--446.

\bibitem[GN]{GrossNebe}
B.H. Gross and G.~Nebe, \emph{Globally maximal arithmetic groups}, J. Algebra
  \textbf{272} (2004), no.~2, 625--642.

\bibitem[Hi]{Hiss}
G.~Hiss, \emph{Die adjungierten {D}arstellungen der {C}hevalley-{G}ruppen},
  Arch. Math. (Basel) \textbf{42} (1984), 408--416.

\bibitem[Ho]{Hogeweij}
G.M.D. Hogeweij, \emph{Almost-classical {L}ie algebras. {I}, {II}}, Nederl.
  Akad. Wetensch. Indag. Math. \textbf{44} (1982), no.~4, 441--460.
  
\bibitem[Hu]{Hum:cc}
J.E. Humphreys, \emph{Conjugacy classes in semisimple algebraic groups}, Math. surveys and monographs \textbf{43}, Amer. Math. Soc., 1995.  

\bibitem[J\,1]{Jantzen}
J.C. Jantzen, \emph{Representations of algebraic groups}, second ed., Math.
  Surveys and Monographs, vol. 107, Amer. Math. Soc., 2003.
  
\bibitem[J\,2]{Jantzen:nil}
\bysame, \emph{Nilpotent orbits in representation theory}, in ``Lie theory: Lie algebras and representations" (J.-P. Anker and B. Orsted, eds.), Progress in Math. \textbf{228}, 2003.

\bibitem[LS]{LaszloSorger}
Y.~Laszlo and C.~Sorger, \emph{The line bundles on the moduli of parabolic
  ${G}$-bundles over curves and their sections}, Ann. Ec. Norm. Sup. (4)
  \textbf{30} (1997), 499--525.
  
\bibitem[Ma]{Mathieu}
O. Mathieu, \emph{Classification des alg\`ebres de Lie simples}, S\'eminaire Bourbaki,  vol. 1998/99, Exp. No. 848, Ast\'erisque, no. 266 (2000), 245--286.

\bibitem[Mer]{MG}
A.S. Merkurjev, \emph{Rost invariants of simply connected algebraic groups}, with a section by S.~Garibaldi,
  in ``Cohomological invariants in Galois cohomology", University Lecture Series,
  vol.~28, Amer.\ Math.\ Soc., 2003.

\bibitem[MPR]{McKPR}
W.G. McKay, J.~Patera, and D.W. Rand, \emph{Tables of representations of simple
  {L}ie algebras. Volume {I}: exceptional simple {L}ie algebras}, Centre de
  Recherches Math{\'e}matiques, Montreal, 1990.
  
\bibitem[P]{Premet:KW}
A.~Premet, \emph{Irreducible representations of Lie algebras of reductive groups and the Kac-Weisfeiler conjecture}, Invent. math. \textbf{121} (1995), 79--117.

\bibitem[R]{Richardson}
R.W. Richardson, \emph{Conjugacy classes in Lie algebras and algebraic groups}, Ann. Math. (2) \textbf{86} (1967), 1--15.

\bibitem[Se]{Sel:mod}
G.B. Seligman, \emph{Modular {L}ie algebras}, Ergebnisse der Mathematik und
  ihrer Grenzgebiete, vol.~40, Springer, 1967.

\bibitem[SpSt]{SpSt}
T.A. Springer and R.~Steinberg, \emph{Conjugacy classes}, Seminar on Algebraic
  Groups and Related Finite Groups (The Institute for Advanced Study,
  Princeton, N.J., 1968/69), Lecture Notes in Math., vol. 131, Springer,
  Berlin, 1970, pp.~167--266.

\bibitem[St\,61]{St:aut}
R.~Steinberg, \emph{Automorphisms of classical {L}ie algebras}, Pacific J.
  Math. \textbf{11} (1961), 1119--1129 [= Collected Papers, pp.~101--111].

\bibitem[St\,62]{St:gen}
\bysame, \emph{G\'en\'erateurs, relations et rev\^{e}tements de groupes
  alg\'ebriques}, Colloque sur la Th\'eorie des Groupes Alg\'ebriques
  (Bruxelles), Centre Belge de Recherches Math., 1962, pp.~113--127 [= Collected Papers, pp.~133--147].

\bibitem[St\,63]{St:rep}
\bysame, \emph{Representations of algebraic groups}, Nagoya Math. J.
  \textbf{22} (1963), 33--56 [= Collected Papers, pp.~149--172].

\bibitem[St\,68]{St}
\bysame, \emph{Lectures on {C}hevalley groups}, Yale University, New Haven,
  Conn., 1968.

\bibitem[St\,75]{St:tor}
\bysame, \emph{Torsion in reductive groups}, Adv. Math. \textbf{15} (1975),
  no.~1, 63--92 [= Collected Papers, pp.~415--444].

\bibitem[Strade]{Strade}
H.~Strade, \emph{Simple {L}ie algebras over fields of positive characteristic.
  {I}}, Expositions in Mathematics, vol.~38, Walter de Gruyter \&
  Co., Berlin, 2004.
  
\bibitem[T]{Ti:Cl}
J.~Tits, \emph{Classification of algebraic semisimple groups}, Algebraic Groups and Discontinuous Subgroups, Proc. Symp. Pure Math., vol. IX, AMS, 1966, pp.~32--62.
  

\end{thebibliography}
